\documentclass{amsart}

\makeatother
\usepackage[maxbibnames=99, backend=biber, style=alphabetic]{biblatex}
\usepackage{amsmath,amsthm,hyperref,amssymb,mathtools,enumitem,tikz-cd}
\usepackage{mathrsfs}
\usepackage{graphicx} 

\title{The derived category of a locally compact space is rarely smooth}
\author{Oscar Harr}
\date{\today}

\newtheorem{thm}{Theorem}
\newtheorem{lem}[thm]{Lemma}

\newtheorem{cor}[thm]{Corollary}

\theoremstyle{definition}
\newtheorem*{defn}{Definition}

\newtheorem{rmk}[thm]{Remark}

\begin{document}
\maketitle
\begin{abstract}
  We show that the derived category of a locally compact Hausdorff space $X$
  is smooth in the sense of non-commutative geometry if and only if $X$ is discrete and finite.
\end{abstract}
The purpose of this short note is to record the following:
\begin{thm}\label{thm:rarely-smooth}
  Let $\mathcal C\in\operatorname{CAlg}\left(\mathcal P\mathrm{r}^{\mathrm{dual}}_{\mathrm{st}}\right)$ be a locally rigid $\infty$-category~\cite{arinkin2022stacklocalsystemsrestricted,ramzi2026locallyrigidinftycategories} such that $\mathcal C^\omega\not\simeq 0$ (e.g.~if $\mathcal C$ is rigid).

  If $X$ is a locally compact Hausdorff space such that the associated $\infty$-category category of $\mathcal C$-valued sheaves $\operatorname{Shv}(X;\mathcal C)$ is smooth, then $X$ is finite.
\end{thm}
Here smoothness is meant in the sense of Kontsevich's non-commutative geometry~\cites[][]{kontsevich-ias}[][Ch~11]{SAG}.
\begin{defn}\label{defn:smooth}
  Let $\mathcal C$ be as in Theorem~\ref{thm:rarely-smooth}. A stable $\mathcal C$-linear $\infty$-category $\mathcal M\in\operatorname{Mod}_{\mathcal C}(\mathcal P\mathrm{r}^{\mathrm{dual}}_{\mathrm{st}})$ is \emph{smooth} (resp. \emph{proper}) if
  \begin{enumerate}[label=(\roman*)]
  \item $\mathcal M$ is dualizable with respect to the $\mathcal C$-linear Lurie tensor product; and
  \item The coevaluation $\mathcal C\to\mathcal M\otimes_{\mathcal C}\mathcal M^\vee$ (resp. the evaluation $\mathcal M^\vee\otimes_{\mathcal C}\mathcal M\to\mathcal C$) is strongly continuous.
  \end{enumerate}
\end{defn}
As a special case of our theorem, we recover
\begin{cor}[Ramzi~\cite{ramzi2022dualizability}]
  If $M$ is a topological manifold such that $\operatorname{Shv}(M;\mathrm{Sp})$ is smooth, then $M$ is discrete.
\end{cor}
We will need the following elementary fact from point-set topology:
\begin{lem}\label{lem:topology-lemma}
  For a topological space $X$, the diagonal $\Delta\subseteq X\times X$ is open if and only if $X$ is discrete.
\end{lem}
\begin{proof}
  If $X$ is discrete then so is $X\times X$, and hence $\Delta\subseteq X\times X$ is open since all subsets of $X\times X$ are open.
  On the other hand, suppose $\Delta$ is open and let $x\in X$ be an arbitrary point. By assumption there is some open neighborhood $W$ of $(x,x)$ with $W\subseteq\Delta$.
  By definition of the product topology, we can find open neighborhoods $U$ and $V$ of $x$ in $X$ such that $U\times V\subseteq\Delta$. But the latter implies $U = V = \lbrace x\rbrace$, finishing the proof.
\end{proof}
\begin{proof}[Proof of Theorem~\ref{thm:rarely-smooth}]
  The coevaluation for $ \operatorname{Shv}(X;\mathcal C)$ is given by the composition
  \[
    \mathcal C\xrightarrow{p^*}\operatorname{Shv}(X;\mathcal C)\xrightarrow{\Delta_*}\operatorname{Shv}(X\times X;\mathcal C)\simeq\operatorname{Shv}(X;\mathcal C)\otimes_{\mathcal C} \operatorname{Shv}(X;\mathcal C),
  \]
  where $p\colon X\to\mathrm{pt}$ is the projection to a point, $\Delta\colon X\hookrightarrow X\times X$ is the diagonal inclusion. Indeed, it follows from e.g.~\cite{volpe2023sixoperations,Mann2022} that the assignments
  \[
    X\mapsto\operatorname{Shv}(X;\mathcal C)
    \quad\text{and}\quad
    \begin{tikzcd}[row sep=small,column sep=tiny]
      & Y\arrow[dl,"f"]\arrow[dr,"f'"] \\
      X && X'
    \end{tikzcd}
    \mapsto (f')_!f^*\colon\operatorname{Shv}(X;\mathcal C)\to \operatorname{Shv}(X';\mathcal C)
  \]
  assemble into a symmetric monoidal six-functor formalism on locally compact Hausdorff spaces, and the claimed formula for the coevaluation is then a general fact about symmetric monoidal six-functor formalisms, see~\cite[p.~xxxvii]{Gaitsgory_Rozenblyum2017}.

  Since $\mathcal C^\omega\not\simeq 0$, we can pick some nonzero compact object $C\in\mathcal C^\omega$. The coevaluation is assumed to be strongly continuous, so it must preserve the compact object $C$. But it takes $C$ to $\Delta_*p^*C$. By proper base change we find for each point $x\colon\mathrm{pt}\to X$ that
  \[
    x^*\Delta_*p^*C\simeq
    \begin{cases}
      C,&\text{if }x\in\Delta,\\
      0,&\text{otherwise.}
    \end{cases}
  \]
  In particular, the support of $\Delta_*p^*C$ is exactly the diagonal $\Delta\subseteq X\times X$. It follows from the characterization of compact objects in $\operatorname{Shv}(X;\mathcal C)$ \cite{Harr2025,Efimov2025} (cf.~the characterization of \emph{internally} compact objects in \cite[Thm~2.5.4.8]{martini2025presentabilitytopoiinternalhigher}) that $\Delta$ must be open and compact in $X\times X$. Lemma~\ref{lem:topology-lemma} now implies that $X$ is discrete. But we also saw that $X\cong\Delta$ is compact, whence $X$ must be finite.
\end{proof}
\begin{rmk}
  Conversely, if $X$ is a finite discrete space and $\mathcal C$ is as in Theorem~\ref{thm:rarely-smooth}, then $\operatorname{Shv}(X;\mathcal C)$ is easily seen to be smooth.
\end{rmk}
\begin{rmk}
  Theorem~\ref{thm:rarely-smooth} stands in stark contrast to the situation in algebraic geometry; the category $\operatorname{QCoh}(X)$ of (derived) quasi-coherent sheaves on a scheme detects smoothness, whereas we have seen that smoothness of $\operatorname{Shv}(X;\mathcal C)$ is a useless condition. On the other hand, $\operatorname{Shv}(X;\mathcal C)$ is proper for a large and important class of spaces~\cite[Prop~3.27]{Efimov2025b}.
\end{rmk}
\subsection*{Acknowledgements}
In an earlier version of this note, I misidentified the coevaluation for $\operatorname{Shv}(X;\mathcal C)$. This led to a longer and more satisfying---but erroneous---proof. I thank the anonymous referree for pointing out this mistake and for other useful comments, and apologize for any confusion caused by my mistake. I am deeply grateful
to Maxime Ramzi for several valuable discussions, remarks and literature recommendations. I was partially supported by the Danish National Research Foundation
through the Copenhagen Centre for Geometry and Topology (DRNF151), and Dan Petersen's Wallenberg Scholar fellowship.
\printbibliography{}
\end{document}